\documentclass[a4paper]{amsart}
\usepackage[utf8]{inputenc}
\usepackage{hyperref}
\usepackage{amsmath,amssymb}
\newcommand{\barMov}{\overline{\textrm{Mov}}}
\newcommand{\Mov}{\textrm{Mov}}
\newtheorem{Theorem}{Theorem}
\newtheorem{Corollary}{Corollary}

\newtheorem{Lemma}{Lemma}
\newtheorem{Conjecture}{Conjecture}

\title[Nonvanishing and Abundance for cones of movable divisors]
{Nonvanishing and Abundance \\ for cones of movable divisors}
\author{G. Bini, M.C. Brambilla, C. Fontanari, E. Postinghel}

\address{Gilberto Bini, Dipartimento di Matematica e Informatica, Universit\`a degli Studi di Palermo, Via Archirafi 34, 90138 Palermo, Italy}
\email{gilberto.bini@unipa.it}

\address{Maria Chiara Brambilla, Universit\`a Politecnica delle Marche, Via Brecce Bianche, 60131 Ancona, Italy}
\email{m.c.brambilla@univpm.it}

\address{Claudio Fontanari, Dipartimento di Matematica, Universit\`a degli Studi di Trento, Via Sommarive 14, 38123 Povo, Trento, Italy}
\email{claudio.fontanari@unitn.it}

\address{Elisa Postinghel, Dipartimento di Matematica, Universit\`a degli Studi di Trento, Via Sommarive 14, 38123 Povo, Trento, Italy}
\email{elisa.postinghel@unitn.it}

 

\begin{document}

\begin{abstract}
\noindent
Let $\barMov^k(X)$ be the closure of the cone $\Mov^k(X)$ generated 
by classes of effective divisors on a projective variety $X$ 
with stable base locus of codimension at least $k+1$. 
We propose a generalized version of the Log Nonvanishing Conjecture
and of the Log Abundance Conjecture for a klt pair $(X,\Delta)$, that is: 
if $K_X+\Delta \in \barMov^{k}(X)$, then $K_X+\Delta \in \Mov^{k}(X)$.  
Moreover, we prove that if the Log Minimal Model Program, the Log Nonvanishing, 
and the Log Abundance hold, then so does our conjecture.
\end{abstract}

\maketitle

\section{Introduction}

Let $X$ be a normal projective variety of dimension $n$ 
and fix an integer $k$ with $0 \le k \le n-1$. 

The cone $\Mov^k(X) \subset N^1(X)_\mathbb{R}$ is the cone generated 
by numerical classes of effective Cartier divisors $D$ such that the stable 
base locus of $\vert D \vert$ has codimension at least $k+1$ in $X$. 
The cone $\barMov^k(X)$ is the closure of $\Mov^k(X)$ in $N^1(X)_\mathbb{R}$. 

In particular, $\Mov^0(X) = \textrm{Eff}(X)$ is the cone of effective divisors 
and $\barMov^0(X) = \overline{\textrm{Eff}}(X)$ is the cone of 
pseudoeffective divisors. Next, $\barMov^1(X)$ is just the standard 
movable cone. Finally, $\Mov^{n-1}(X) = \textrm{Sem}(X)$ 
turn out to be the cone of semiample divisors by the Theorem of Zariski-Fujita (see 
for instance \cite{Laz}, Remark 2.1.32) and $\barMov^{n-1}(X) = {\textrm{Nef}}(X)$ 
is the nef cone. 

The cones $\Mov^k(X)$ form a natural filtration: 
$$
\textrm{Sem}(X) = \Mov^{n-1}(X) \subseteq \Mov^{n-2}(X) \subseteq \ldots \subseteq
\Mov^1(X) \subseteq \Mov^0(X) = \textrm{Eff}(X)
$$
Taking closures, we obtain a similar filtration of the corresponding closed cones. 

In \cite{bdps}, the cones $\barMov^k(X)$ were denoted by $\mathcal{D}_k(X)$ and 
called \emph{the cones of divisors that are ample in codimension $k$}, following the 
terminology introduced in \cite{Payne}. 

If $X$ is a Mori Dream Space, then we have $\Mov^k(X) = \barMov^k(X)$ for $k \in \{0,1,n-1 \}$ by definition and the same holds also for all other values of $k$ by \cite{bdps}, Theorem 1.1. For an arbitrary $X$, the cones $\Mov^k(X)$ need not be closed, but the canonical divisor $K_X$ is expected to behave well at least with respect to $\Mov^0(X)$ and $\Mov^{n-1}(X)$. Indeed, we have the following celebrated conjectures in birational geometry: 

\begin{Conjecture} \emph{(Log Nonvanishing)}
Let $(X, \Delta)$ be a projective klt pair.   
If $K_X+\Delta \in \overline{\textrm{Eff}}(X) = \barMov^0(X)$,
then $K_X+\Delta \in \textrm{Eff}(X) = \Mov^0(X)$.
\end{Conjecture}

\begin{Conjecture} \emph{(Log Abundance)}
Let $(X, \Delta)$ be a projective klt pair with $\dim(X)=n$.   
If $K_X+\Delta \in {\textrm{Nef}}(X) = \barMov^{n-1}(X)$, 
then $K_X+\Delta \in \textrm{Sem}(X) = \Mov^{n-1}(X)$.
\end{Conjecture}

Both conjectures are true for $n \le 3$ by \cite{KMM} and 
for $K_X+\Delta$ big and arbitrary $n$ by Kawamata-Shokurov 
Basepoint free theorem, but are widely open in general (see 
for instance \cite{LP}). 

In this paper we introduce a natural generalization of the above conjectures 
for $\barMov^{k}(X)$ with $0 \le k \le n-1$. Namely, we propose the following:

\begin{Conjecture} \label{generalized}
Let $(X, \Delta)$ be a projective klt pair with $\dim(X)=n$ 
and fix an integer $k$ with $0 \le k \le n-1$.    
If $K_X+\Delta \in \barMov^{k}(X)$, then $K_X+\Delta \in \Mov^{k}(X)$.
\end{Conjecture}

The statement of Conjecture \ref{generalized} for $1 \le k \le n-2$ is 
formally the same as the one for $k=0$ and for $k=n-1$, but it turns 
out to be substantially weaker. Indeed, the assumption in the case 
$k=n-1$ is very explicit: $K_X+\Delta \in {\textrm{Nef}}(X) = \barMov^{n-1}(X)$
means that the degree of $K_X+\Delta$ on every irreducible curve on $X$ 
is nonnegative, and the same is true also in the case $k=0$: indeed, 
$K_X+\Delta \in \overline{\textrm{Eff}}(X) = \barMov^0(X)$ means 
precisely that the degree of $K_X+\Delta$ on any irreducible member of a family 
of curves covering $X$ is nonnegative, by the main result in \cite{BDPP}. 
On the other hand, such a \emph{strong duality} between divisors and 
curves on $X$ may fail in the intermediate cases $1 \le k \le n-1$ even 
for a toric threefold (see \cite{Payne}, Example 1). A weaker form 
of duality, involving divisors and curves on the (small) birational models of $X$, 
still holds for any $k$, at least for Fano type varieties and Mori Dream Spaces (see \cite{Choi}, Theorem 1.1 and Corollary 1.1), and more generally 
for smooth projective varieties (see  
\cite{Lehmann}, Corollary 3.3 and Remark 3.4), 
but such a formulation is 
admittedly less effective than in the two classical cases $k=n-1$ and $k=0$. 

Our main result concerning Conjecture \ref{generalized} is the following:

\begin{Theorem}\label{main}
Let $(X, \Delta)$ be a projective klt pair with $\dim(X)=n$ 
and fix an integer $k$ with $0 \le k \le n-1$.    
Assume that the Log Minimal Model Program, the Log Nonvanishing Conjecture, 
and the Log Abundance Conjecture hold for 
every pair $(X', \Delta')$ with 
$f: X \dashrightarrow X'$ an isomorphism in codimension $k$ and $\Delta' = f_*\Delta$.
If $K_X+\Delta \in \barMov^{k}(X)$, then $K_X+\Delta \in \Mov^{k}(X)$.
\end{Theorem}

As a consequence of Mori theory in dimension three, we obtain the following unconditional
result: 

\begin{Corollary}
Let $(X, \Delta)$ be a projective klt pair with $\dim(X)=3$.    
If $K_X+\Delta \in \barMov^{1}(X)$, 
then $K_X+\Delta \in \Mov^{1}(X)$.
\end{Corollary}

In particular, by taking $K_X = 0$ and by replacing an effective divisor $D$ with 
$\varepsilon D$ with $0 < \varepsilon << 1$, we deduce for K-trivial threefolds 
that an effective divisor lying on the boundary of $\barMov^{1}(X)$ is indeed movable, 
exactly as in the case of Fano varieties (or, more generally, of Mori Dream Spaces):

\begin{Corollary}
Let $X$ be a projective normal threefold with trivial canonical bundle.    
Then we have $\barMov^{1}(X) \cap \textrm{Eff}(X) = \Mov^{1}(X)$.
\end{Corollary}

More generally, as a consequence of Theorem \ref{main}, the same equality $\barMov^{1}(X) \cap \textrm{Eff}(X) = \Mov^{1}(X)$ holds for projective normal varieties 
$X$ of arbitrary dimension with trivial canonical bundle, provided that the Log Minimal Model Program, the Log Nonvanishing Conjecture, and the Log Abundance Conjecture hold for every small birational modification $X'$ of $X$. This fact is remarkable because in the statement of the Morrison-Kawamata Cone Conjecture (see for instance \cite{LOP}, Conjecture 3.3 (2)) the cone $\barMov^{1}(X) \cap \textrm{Eff}(X)$ plays a crucial role. 

Finally, since by \cite{BCHM} the Log Minimal Model Program holds for varieties of log general type, we get: 

\begin{Corollary}
Let $(X, \Delta)$ be a projective klt pair with $\dim(X)=n$ 
and fix an integer $k$ with $0 \le k \le n-1$.    
If $K_X+\Delta \in \barMov^{k}(X)$ is big, then 
$K_X+\Delta \in \Mov^{k}(X)$.
\end{Corollary}

We work over the complex field $\mathbb{C}$.

\section{Proof of Theorem \ref{main}}
Since $K_X+\Delta \in \barMov^{k}(X) \subseteq \barMov^{0}(X) = \overline{\textrm{Eff}}(X)$, 
by the Log Nonvanishing assumption on $(X, \Delta)$ we have $K_X+\Delta \in \textrm{Eff}(X)$. 

We claim that there exists $f: X \dashrightarrow  X'$ an isomorphism 
in codimension $k$ (i.e. a birational map with exceptional locus of 
codimension at least $k+1$ in $X$) with $K_{X'} = f_*K_X$
such that $K_{X'}+ \Delta'$ is nef, where $\Delta' = f_*\Delta$. 
Indeed, by the Log Minimal Model Program assumption on $(X, \Delta)$ 
there exists a birational map $f: X \dashrightarrow X'$ with 
$K_{X'} = f_*K_X$ and $K_{X'}+ \Delta'$ nef. 
If by contradiction there is a subvariety $V$ in the exceptional 
locus of $f$ of codimension at most $k$ in $X$, then there is a 
$k$-moving curve class $C$ sweeping out $V$ such that $(K_X+\Delta).C < 0$.
Now, since $K_X+\Delta \in \barMov^{k}(X)$, there exists a sequence of 
$\mathbb{R}$-divisors $D_i \in \Mov^{k}(X)$ converging to $K_X+\Delta$; 
in particular, we have $D_i.C < 0$ for any $i$ sufficiently large. 
It follows that $V$ is contained in the (stable) base locus of 
$D_i$, a contradiction.

Next, by the Log Abundance assumption on $(X', \Delta')$, since
$f_*(K_X+ \Delta) = K_{X'}+ \Delta'$ is nef then it is semiample, 
in particular $f_*(K_X+ \Delta) \in \Mov^{k}(X')$. Hence we may 
conclude by invoking the previous claim and applying the following 
Lemma to $D := K_X+ \Delta$.

\begin{Lemma}
Let $D$ be an effective Cartier divisor on a normal projective 
variety $X$ and let $f: X \dashrightarrow  X'$ be an isomorphism 
in codimension $k$. If $f_*D \in \Mov^{k}(X')$, then $D \in \Mov^{k}(X)$.
\end{Lemma}

\proof
By contradiction, assume that $D \notin \Mov^{k}(X)$, so that 
there exists a subvariety $W$ in the stable base locus of $D$ 
of codimension at most $k$ in $X$. Since by assumption the 
exceptional locus of $f$ has codimension at least $k+1$ 
in $X$, it follows that $f_{\vert W}: W \dashrightarrow f(W)$
is birational, in particular we have $\dim f(W) = \dim W$. 
We claim that $f(W)$ is contained in the stable base locus 
of $f_*(D)$. Indeed, we have $\vert f_*(D) \vert = f_* \vert D \vert$
since $f$ is an isomorphism in codimension $k \ge 1$ and $f(W)$
is contained in every element of $f_* \vert D \vert$ since 
$W$ is contained in every element of $\vert D \vert$.
From the claim it follows that $f_*D \notin \Mov^{k}(X')$, 
a contradiction.

\qed

\section*{Acknowledgements}
We are grateful to Paolo Cascini and Vladimir Lazi\'c for sharing their realiable insights on birational geometry with us.
The authors are members of GNSAGA of INdAM (Italy). This project has been
developed in April 2024 as a Research in Pairs program at CIRM, Trento (Italy). 
The forth author's research is funded by the European Union under the Next Generation EU PRIN 2022 \emph{Birational geometry of moduli spaces and special varieties}, Prot. n. 20223B5S8L.

\end{document}